\documentclass[12]{article}
\usepackage{bm}
\usepackage{fleqn}
\usepackage{graphicx}
\begin{document}
\Large
\begin{center}
\bf{Combinatorial Intricacies of Labeled Fano Planes}
\end{center}
\vspace*{-.3cm}
\begin{center}
Metod Saniga
\end{center}
\vspace*{-.5cm} \normalsize
\begin{center}
Astronomical Institute, Slovak Academy of Sciences\\
SK-05960 Tatransk\' a Lomnica, Slovak Republic\\
(msaniga@astro.sk)  
\end{center}

\vspace*{-.4cm} \noindent \hrulefill

\vspace*{-.1cm} \noindent {\bf Abstract}

\noindent
Given a seven-element set $X = \{1,2,3,4,5,6,7\}$, there are 30 ways to define a Fano plane on it. Let us call a line of such Fano plane, that is to say an unordered triple from $X$, ordinary or defective according as the sum of two smaller integers from the triple is or is not equal to the remaining one, respectively. A point of the labeled Fano plane is said to be of order $s$, $0 \leq s \leq 3$, if there are $s$ {\it defective} lines passing through it. With such structural refinement in mind, the 30 Fano planes are shown to fall into eight distinct types. Out of the total of 35 lines, nine ordinary lines are of five different kinds, whereas the remaining 26 defective lines yield as many as ten distinct types. It is shown, in particular, that no labeled Fano plane can have all points of zeroth order, or feature just one point of order two. A connection with prominent configurations in Steiner triple systems is also pointed out.\\

\vspace*{-.2cm}


\vspace*{-.2cm} \noindent \hrulefill

\vspace*{.1cm}
\noindent
It is well known \cite{hall,lr,pol1} that there are thirty different Fano planes on a given seven-element set. Slightly rephrased, there are thirty different ways to label the points of the Fano plane by integers from 1 to 7, two labeled Fano planes having zero, one or three lines in common and each line occurring in six Fano planes.
The set of thirty labeled Fano planes can be uniquely partitioned into two sets, say $A$ and $B$, of fifteen elements each, such that any two labeled Fanos in the same set have just one line in common.  Following the labeling adopted by Polster \cite{pol2}, the two sets read 
$$~A_1: \{124,136,157,\overline{235},267,\overline{347},456\},\hspace*{1.0cm}B_1: \{127,136,\overline{145},\overline{235},\overline{246},\overline{347},567\},$$
$$~A_2:~ \{127,136,\overline{145},234,256,357,467\},\hspace*{1.0cm}B_2: \{125,136,147,234,267,357,456\},$$
$$~A_3:~ \{125,136,147,237,\overline{246},345,567\},\hspace*{1.0cm}B_3: \{124,136,157,237,256,345,467\},$$
$$~A_4:~ \{125,\overline{134},\overline{167},236,247,357,456\},\hspace*{1.0cm}B_4: \{127,\overline{134},\overline{156},236,245,357,467\},$$
$$~A_5:~ \{127,135,146,236,245,\overline{347},567\},\hspace*{1.0cm}B_5: \{124,135,\overline{167},236,\overline{257},\overline{347},456\},$$
$$~A_6:~ \{124,137,\overline{156},236,\overline{257},345,467\},\hspace*{1.0cm}B_6: \{125,137,146,236,247,345,567\},$$
$$~A_7:~ \{\overline{123},147,\overline{156},245,267,346,357\},\hspace*{1.0cm}B_7: \{\overline{123},\overline{145},\overline{167},247,256,346,357\},$$
$$~A_8:~ \{124,135,\overline{167},237,256,346,457\},\hspace*{1.0cm}B_8: \{126,135,147,237,245,346,567\},$$
$$~A_9:~ \{126,137,\overline{145},\overline{235},247,346,567\},\hspace*{1.0cm}B_9: \{124,137,\overline{156},\overline{235},267,346,457\},$$
$$A_{10}:~ \{\overline{123},\overline{145},\overline{167},\overline{246},\overline{257},\overline{347},356\},\hspace*{1.0cm}B_{10}: \{\overline{123},146,157,245,267,\overline{347},356\},$$
$$A_{11}:~ \{126,\overline{134},157,237,245,356,467\},\hspace*{1.0cm}B_{11}: \{125,\overline{134},\overline{167},237,\overline{246},356,457\},$$
$$A_{12}:~ \{125,137,146,234,267,356,457\},\hspace*{1.0cm}B_{12}: \{126,137,\overline{145},234,\overline{257},356,467\},$$
$$A_{13}:~ \{\overline{123},146,157,247,256,345,367\},\hspace*{1.0cm}B_{13}: \{\overline{123},147,\overline{156},\overline{246},\overline{257},345,367\},$$
$$A_{14}:~ \{127,\overline{134},\overline{156},\overline{235},\overline{246},367,457\},\hspace*{1.0cm}B_{14}: \{126,\overline{134},157,\overline{235},247,367,456\},$$
$$A_{15}:~ \{126,135,147,234,\overline{257},367,456\},\hspace*{1.0cm}B_{15}: \{127,135,146,234,256,367,457\},$$
\noindent
where for triples of integers we use a shorthand notation $\{a,b,c\} = abc$.

Given a line $abc$, where, without loss of generality, we can take $1 \leq  a < b < c \leq 7$, we shall distinguish between the cases when $a+b = c$ and $a+b \neq c$ and call the former/latter ordinary/defective \cite{cd}; for the sake of convenience, in the above-given sets $A$ and $B$ all ordinary lines are denoted by overbars. Furthermore, a point of the labeled Fano plane is said to be of order $s$, $0 \leq s \leq 3$, if there are $s$ {\it defective} lines passing through it; hence, in addition to two different kinds of lines, a labeled Fano plane can potentially feature up to four distinct types of points. A detailed inspection of each of the 30 labeled Fano planes given above shows that they fall, in terms of such structural refinement, into eight different types as  portrayed in Figure 1 and summarized in Table 1. It is also interesting to compare the distributions of types within each 15-element set, given in Tables 2 and 3, as well as cardinalities of individual types, listed in Table 4; note a pronounced asymmetry between sets $A$ and $B$.

\begin{figure}[pth!]
\centerline{\includegraphics[width=11.2truecm,clip=]{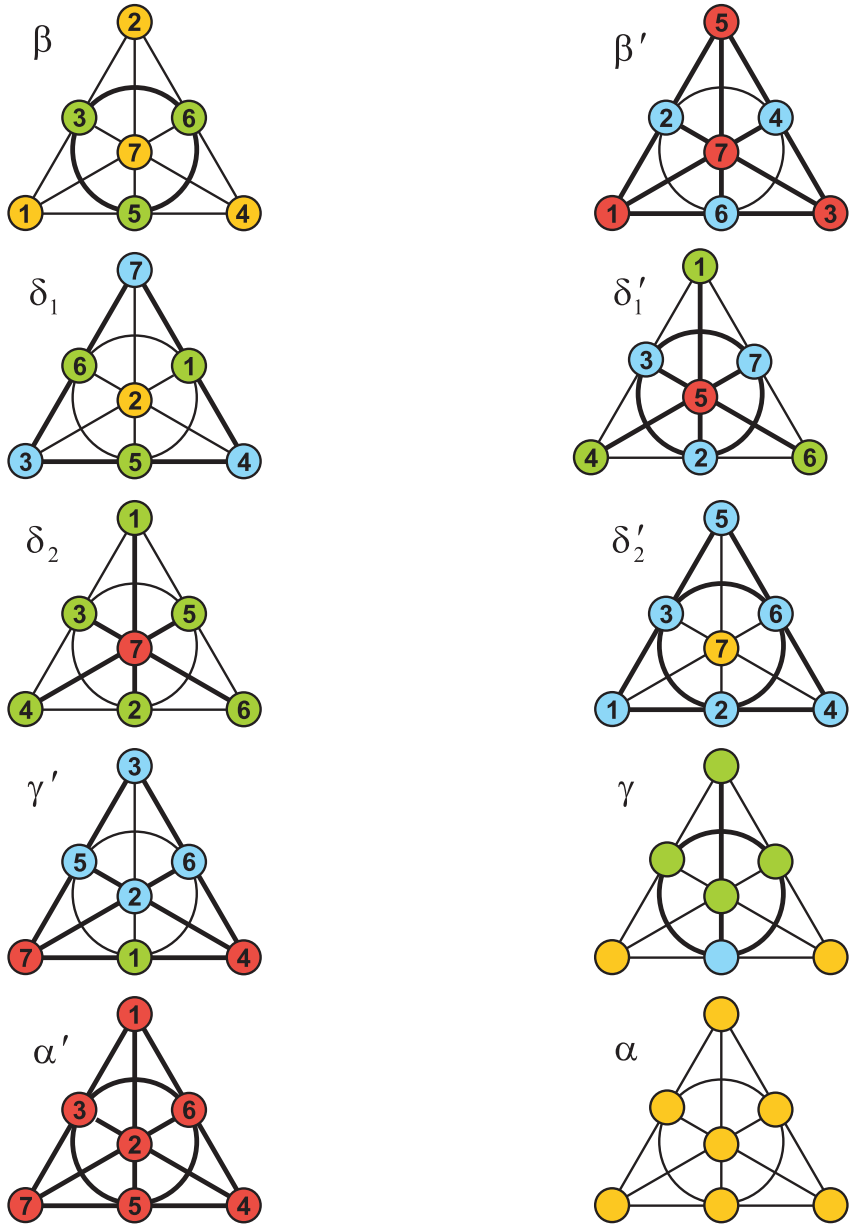}}
\caption{A diagrammatic illustration of representatives of eight distinct types of numbered Fano planes. A point of order three, two, one and/or zero is represented by red, blue, green and/or yellow color, respectively; heavy lines are defective. The types $\alpha$ and $\gamma$ do not exist, i.\,e. it is impossible label the Fano plane in such a way that all or just five of its lines would be ordinary.}
\end{figure}

\begin{table}[t]
\begin{center}
\caption{The eight distinct types of labeled Fano planes, each being uniquely characterized by the number of points of every particular order.  To understand better the symmetry pattern, we arrange complementary types into pairs and list also the two non-existent types $\alpha$ and $\gamma$.} 
\vspace*{0.2cm}
{\begin{tabular}{|c|c|c|c|c|} \hline \hline
\multicolumn{1}{|c|}{} &  \multicolumn{4}{|c|}{}   \\
\multicolumn{1}{|c|}{} &   \multicolumn{4}{|c|}{Points of Order} \\
 \cline{2-5}
Type         &   0    &   1    &   2   &  3   \\
\hline
($\alpha$)   &  (7)   &  (0)   &  (0)  & (0)    \\
$\alpha'$    &   0    &   0    &   0   &  7    \\
\hline
$\beta$      &   4    &   3    &   0   &  0    \\
$\beta'$     &   0    &   0    &   3   &  4    \\
\hline
($\gamma$)   &  (2)   &  (4)   &  (1)  & (0)    \\
$\gamma'$    &   0    &   1    &   4   &  2    \\
\hline
$\delta_1$   &   1    &   3    &   3   &  0    \\
$\delta_1'$  &   0    &   3    &   3   &  1    \\
\hline
$\delta_2$   &   0    &   6    &   0   &  1    \\
$\delta_2'$  &   1    &   0    &   6   &  0    \\
\hline \hline
\end{tabular}}
\end{center}
\end{table}

\begin{table}[pth!]
\begin{center}
\caption{Types of Fano planes in set $A$. } \vspace*{0.4cm}
\begin{tabular}{|c|ccccccccccccccc|}
\hline \hline
 Plane & 1 & 2 & 3 & 4 & 5 & 6 & 7 & 8 & 9 & 10 & 11 & 12 & 13 & 14 & 15  \\
\hline
 Type & $\gamma'$ & $\beta'$ & $\beta'$ & $\gamma'$ & $\beta'$ & $\gamma'$ & $\gamma'$ & $\beta'$ & $\gamma'$ & $\beta$ & $\beta'$ & $\alpha'$ & $\beta'$ & $\delta_2$ & $\beta'$  \\
 \hline \hline
\end{tabular}
\end{center}
\end{table}

\begin{table}[pth!]
\begin{center}
\caption{Types of Fano planes in set $B$. } \vspace*{0.4cm}
\begin{tabular}{|c|ccccccccccccccc|}
\hline \hline
 Plane & 1 & 2 & 3 & 4 & 5 & 6 & 7 & 8 & 9 & 10 & 11 & 12 & 13 & 14 & 15  \\
\hline
 Type    & $\delta_1$  & $\alpha'$  & $\alpha'$  & $\gamma'$ & $\delta_2'$  & $\alpha'$  & $\delta_2'$  & $\alpha'$  & $\gamma'$ & $\gamma'$ & $\delta_1'$  & $\gamma'$ & $\delta_1$  & $\gamma'$ & $\alpha'$   \\
 \hline \hline
\end{tabular}
\end{center}
\end{table}
 
\begin{table}[pth!]
\begin{center}
\caption{Cardinalities of individual types of Fano planes. } \vspace*{0.4cm}
\begin{tabular}{|l|cccccccccc|}
\hline \hline
  Type & ($\alpha$) & $\beta$  & ($\gamma$)  & $\delta_1$ & $\delta_2$ & $\alpha'$ & $\beta'$  & $\gamma'$  & $\delta_1'$ & $\delta_2'$  \\
\hline
Set $A$   & $-$ & 1 & $-$ & 0 & 1 & 1 & 7 &  5 & 0 & 0 \\
Set $B$   & $-$ & 0 & $-$ & 2 & 0 & 5 & 0 &  5 & 1 & 2 \\
\hline
Total & $-$ & 1 & $-$ & 2 & 1 & 6 & 7 & 10 & 1 & 2 \\
 \hline \hline
\end{tabular}
\end{center}
\end{table}
In addition to classifying planes, we can also classify lines in terms of the types of six planes passing through each of them, the corresponding findings being given in Table 5. We find that nine ordinary lines are of five different kinds, whereas the remaining 26 defective lines yield as many as ten distinct types. We further observe that no ordinary line exhibits a plane of type $\alpha'$ and there is only one defective line that features plane of type $\beta$. It is also worth pointing out that it is only one plane in set $A$ that is devoid of ordinary lines, whereas in set $B$ one finds four such planes.
\begin{figure}[pth!]
\centerline{\includegraphics[width=14truecm,clip=]{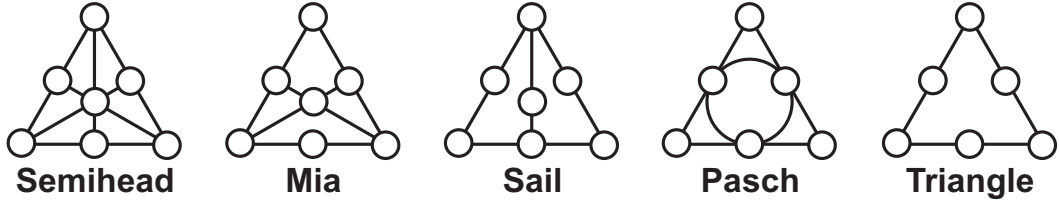}}
\caption{Five distinguished Fano derivatives, with their traditional names, found within labeled Fano planes.}
\end{figure}
\begin{table}[pth!]
\begin{center}
\caption{The types of lines; nine ordinary lines go first. Lines belong to a given type if they possess  the same string of parameters.} \vspace*{0.4cm}
\begin{tabular}{|l|cccccccccc|}
\hline \hline
  Line & ($\alpha$) & $\beta$  & ($\gamma$)  & $\delta_1$ & $\delta_2$ & $\alpha'$ & $\beta'$  & $\gamma'$  & $\delta_1'$ & $\delta_2'$  \\
\hline
123   & $-$ & 1 & $-$ & 1 & 0 & 0 & 1 &  2 & 0 & 1 \\
145   & $-$ & 1 & $-$ & 1 & 0 & 0 & 1 &  2 & 0 & 1 \\
257   & $-$ & 1 & $-$ & 1 & 0 & 0 & 1 &  2 & 0 & 1 \\
347   & $-$ & 1 & $-$ & 1 & 0 & 0 & 1 &  2 & 0 & 1 \\
\hline
156   & $-$ & 0 & $-$ & 1 & 1 & 0 & 0 &  4 & 0 & 0 \\
235   & $-$ & 0 & $-$ & 1 & 1 & 0 & 0 &  4 & 0 & 0 \\
\hline
246   & $-$ & 1 & $-$ & 2 & 1 & 0 & 1 &  0 & 1 & 0 \\
\hline
167   & $-$ & 1 & $-$ & 0 & 0 & 0 & 1 &  1 & 1 & 2 \\
\hline
134   & $-$ & 0 & $-$ & 0 & 1 & 0 & 1 &  3 & 1 & 0 \\
\hline \hline
124                & $-$ & 0 & $-$ & 0 & 0 & 1 & 1 &  3 & 0 & 1 \\
236                & $-$ & 0 & $-$ & 0 & 0 & 1 & 1 &  3 & 0 & 1 \\
247                & $-$ & 0 & $-$ & 0 & 0 & 1 & 1 &  3 & 0 & 1 \\
346                & $-$ & 0 & $-$ & 0 & 0 & 1 & 1 &  3 & 0 & 1 \\
357                & $-$ & 0 & $-$ & 0 & 0 & 1 & 1 &  3 & 0 & 1 \\
456                & $-$ & 0 & $-$ & 0 & 0 & 1 & 1 &  3 & 0 & 1 \\
\hline
136                & $-$ & 0 & $-$ & 1 & 0 & 2 & 2 &  1 & 0 & 0 \\
147                & $-$ & 0 & $-$ & 1 & 0 & 2 & 2 &  1 & 0 & 0 \\
345                & $-$ & 0 & $-$ & 1 & 0 & 2 & 2 &  1 & 0 & 0 \\
567                & $-$ & 0 & $-$ & 1 & 0 & 2 & 2 &  1 & 0 & 0 \\
\hline
126                & $-$ & 0 & $-$ & 0 & 0 & 1 & 2 &  3 & 0 & 0 \\
157                & $-$ & 0 & $-$ & 0 & 0 & 1 & 2 &  3 & 0 & 0 \\
245                & $-$ & 0 & $-$ & 0 & 0 & 1 & 2 &  3 & 0 & 0 \\
467                & $-$ & 0 & $-$ & 0 & 0 & 1 & 2 &  3 & 0 & 0 \\
\hline
135                & $-$ & 0 & $-$ & 0 & 0 & 2 & 3 &  0 & 0 & 1 \\
237                & $-$ & 0 & $-$ & 0 & 0 & 2 & 3 &  0 & 0 & 1 \\
256                & $-$ & 0 & $-$ & 0 & 0 & 2 & 3 &  0 & 0 & 1 \\
\hline
127                & $-$ & 0 & $-$ & 1 & 1 & 1 & 2 &  1 & 0 & 0 \\
367                & $-$ & 0 & $-$ & 1 & 1 & 1 & 2 &  1 & 0 & 0 \\
\hline
146                & $-$ & 0 & $-$ & 0 & 0 & 3 & 2 &  1 & 0 & 0 \\
234                & $-$ & 0 & $-$ & 0 & 0 & 3 & 2 &  1 & 0 & 0 \\
\hline
137                & $-$ & 0 & $-$ & 0 & 0 & 2 & 0 &  4 & 0 & 0 \\
267                & $-$ & 0 & $-$ & 0 & 0 & 2 & 0 &  4 & 0 & 0 \\
\hline
356                & $-$ & 1 & $-$ & 0 & 0 & 1 & 1 &  2 & 1 & 0 \\
\hline
457                & $-$ & 0 & $-$ & 0 & 1 & 2 & 1 &  1 & 1 & 0 \\
\hline
125                & $-$ & 0 & $-$ & 0 & 0 & 3 & 1 &  1 & 1 & 0 \\
\hline \hline
\end{tabular}
\end{center}
\end{table}

To conclude the paper, let us have a look at Figure 1 and focus on that sub-geometry in each labeled Fano plane representative which is formed by defective lines and points lying on them. We shall find the Fano plane itself and, using the language of \cite{ggr}, all its derivatives depicted in Figure 2 that play a crucial role in classifying Steiner triple systems (see, for example, \cite{hor,fgg}); as readily discerned from Figure 1, the semihead lies in a type-$\beta'$ plane, the mia in a type-$\gamma'$ plane, a sail is hosted by the type-$\delta_1'$ plane, the Pasch configuration is located in a type-$\delta_2'$ plane, the triangle in a type-$\delta_1$ plane and the Fano plane itself coincides with a type-$\alpha'$ plane. It should also be pointed out that if one instead considers sub-geometries formed by ordinary lines, then both the Fano plane and its derivate mia are absent.

Obviously, any finite geometry (point-line incidence structure) with lines of size three can be looked at this way. We have already performed the corresponding examination of the M\"obius-Kantor $(8_3)$-configuration, the Pappus $(9_3)$-configuration as well as the Desargues $(10_3)$-configuration and plan to publish its outcome in a separate paper. 

\section*{Acknowledgment}
This work was partially supported by the VEGA Grant Agency, Project 2/0003/13, as well as by the Austrian Science Fund (Fonds zur F\"orderung der Wissenschaftlichen Forschung (FWF)), Research Project M1564--N27. I am grateful to my friend Petr Pracna for electronic versions of the figures.

\end{document}